\documentclass[12pt]{article}
\usepackage{amsxtra,amssymb,amsthm,amsmath,latexsym}

\theoremstyle{plain}

\def\R{{\mathbb R}}

\def\oH{{\overset{\circ}{H}}}
\def\oH1{{\overset{\circ}{H}\kern-.02in{}^1}}

\def\bee{\begin{equation*}}
\def\eee{\end{equation*}}
\def\be{\begin{equation}}
\def\ee{\end{equation}}

\begin{document}

\title{Perturbation of zero surfaces
}

\author{Alexander G. Ramm\\
 Department  of Mathematics, Kansas State University, \\
 Manhattan, KS 66506, USA\\
ramm@math.ksu.edu\\
http://www.math.ksu.edu/\,$\widetilde{\ }$\,ramm
}

\date{}
\maketitle\thispagestyle{empty}

\begin{abstract}

It is proved that if a smooth function $u(x)$, $x\in \R^3$, such that $\inf_{s\in S}|u_N(s)|>0$, where $u_N$
is the normal derivative of $u$ on $S$, has a closed smooth surface $S$ of zeros, then
the function $u(x)+\epsilon v(x)$ has also a closed smooth surface $S_\epsilon$ of zeros. Here $v$ is a smooth function
and $\epsilon>0$ is a sufficiently small number.

\end{abstract}

\section{Introduction}\label{S:1}
 Let $D\subset \mathbb{R}^3$ be a bounded domain containing inside a connected closed $C^3-$smooth surface $S$,
which is the set of zeros of a function $u\in C^3(D)$, so that
\be\label{e1}
 u|_{S}=0.
\ee
Let $N=N_s$ be the unit normal to $S$, such that $u_N=|\nabla u(s)|$, where $u_N$ is the normal derivative of $u$ on $S$.
Let $u_\epsilon:=u+\epsilon v$, where $v\in C^3(D)$ and $\epsilon>0$ is sufficiently small. Assume that
\be\label{e2}
\inf_{s\in S}|\nabla u(s)|\ge 2c_1>0, \qquad c_1=const>0.
\ee
The purpose of this paper is to prove Theorem 1.

{\bf Theorem 1.}  {\em Under the above assumptions there exists a smooth closed surface $S_\epsilon$ such that
$u_\epsilon=0$ on $S_\epsilon$.
 }

In Section 2   Theorem 1 is proved.

Although there are many various results on perturbation theory, see \cite{K}, \cite{R}, the result formulated in Theorem 1 is new.

\section{ Proof of Theorem 1}\label{S:2}
Consider the following equation for $t$:
\be\label{e3}
u(s+tN)+\epsilon v(s+tN)=0,
\ee
where $N=N(s)$ is the normal to $S$ at the point $s$ and $t$ is a parameter. Using the Taylor's formula and relation  \eqref{e1}, one gets from \eqref{e3}
\be\label{e4}
t \Big(\nabla u(s)\cdot N+\epsilon \nabla v(s)\cdot N\Big) +\epsilon v(s) +t^2 \phi=0,
\ee
where $t^2\phi$ is the Lagrange remainder in the Taylor's formula and
\be\label{e5}
 \phi:= \sum_{i,j=1}^3 [u_{x_i x_j}(s+\theta t N)+\epsilon v_{x_i x_j}(s+\theta t N)]N_i N_j, \quad \theta\in (0,1).
\ee
Since the functions $u$ and $v$ belong to $C^3(D)$, the function $\phi=\phi(t,s,\epsilon)$ has a bounded derivative with respect to $t$ uniformly with respect to $s\in S$ and $\epsilon\in (0, 1]$.

Consider equation \eqref{e4} as an equation for $t=t(s)$ in the space $C(S)$. Rewrite \eqref{e4} as
\be\label{e6}
t= -\epsilon \Big(\nabla u(s)\cdot N+\epsilon \nabla v(s)\cdot N \Big)^{-1} v(s) -t^2 \phi \Big(\nabla u(s)\cdot N+\epsilon \nabla v(s)\cdot N \Big)^{-1}:=Bt.
\ee
Let us check that the operator $B$ satisfies the contraction mapping theorem in the set
\be\label{e7}
M:=\{t: max_{s\in S}|t(s)-\epsilon \Big(\nabla u(s)\cdot N+\epsilon \nabla v(s)\cdot N \Big)^{-1} v(s)|\le \delta\},
\ee
where $\delta>0$ is a small number, and $M\in C(S)$.

First, one should check that $B$ maps $M$ into itself. One has
\be\label{e8}
 max_{s\in S}|B t(s)-\epsilon \Big(\nabla u(s)\cdot N+\epsilon \nabla v(s)\cdot N \Big)^{-1} v(s)|\le  max_{s\in S} \frac {t^2 |\phi|}
 {\nabla u(s)\cdot N+\epsilon \nabla v(s)\cdot N}.
\ee
We have chosen $N$ so that $\nabla u(s)\cdot N=|\nabla u(s)|$. This is possible because equation \eqref{e1} implies that
$\nabla u(s)$ is orthogonal to $S$ at the point $s\in S$. Assumption \eqref{e2} implies that for sufficiently small
$\epsilon$ one has
\be\label{e9}
\inf_{s\in S}|\nabla u_\epsilon(s)|\ge c_1.
\ee
Since $\phi$ is continuously differentiable, one has
\be\label{e10}
\sup_{s\in S, t\in (0,1)}|\phi(t,s, \epsilon)|\le c_2.
\ee
Therefore, if
\be\label{e11}
|t(s)|\le \delta,
\ee
then
\be\label{e12}
\frac {t^2(s)|\phi(t,s,\epsilon)|}{|\nabla u(s)|+\epsilon \nabla v(s)\cdot N}\le \frac {c_2}{c_1}\delta^2\le \delta,
\ee
provided that
\be\label{e13}
\frac {c_2}{c_1} \delta\le 1.
\ee
Thus, if \eqref{e13} holds then $B$ maps $M$ into itself.

Let us check that $B$ is a contraction mapping on $M$. One has
\be\label{e14}
|Bt_1-Bt_2|\le c_1^{-1}|t_1^2\phi(t_1,s, \epsilon)-t_2^2\phi(t_2,s,\epsilon)|\le c_3|t_1-t_2|,
\ee
where $c_3\in (0,1)$ if $\delta$ is sufficiently small.
Indeed,
\be\label{e15}
c_3= max_{s\in S, t\le \delta} \Big(2t|\phi(t,s,\epsilon)|+t^2 |\frac {\partial \phi}{\partial t}| \Big)\le c_4\delta<1,
\ee
if $\delta$ is sufficiently small. Here $c_4$ is a constant.

Thus, $B$ is a contraction on $M$. By the contraction mapping principle, equation \eqref{e6} is uniquely
solvable for $t$. Its solution $t=t(s)$ allows one to construct the zero surface $S_\epsilon$ of the function $u_\epsilon$ by the equation  $r=s+t(s)N$, where $r=r(s)$ is the radius vector of the points on $S_\epsilon$.

Theorem 1 is proved. \hfill$\Box$

{\bf Remark 1.} Condition \eqref{e2} is a sufficient condition for the validity of Theorem 1. Although this condition is not necessary, if it does not hold one can construct counterexamples to the conclusion of Theorem 1. For example, assume that $u(x)\ge 0$ in $\R^3$
and $u(x)=0$ on $S$,  let $v> 0$ in $\R^3$ and $\epsilon>0$. Then the function $u_\epsilon=u+\epsilon v$ does not have
zeros in $\R^3$.

{\bf Remark 2.} In scattering theory the following question is of interest: assume that $u(x)$ is an entire function of exponential type, $u(x)=\int_{S^2} e^{ik\beta \cdot x}f(\beta)d\beta$, where $f\in L^2(S^2)$, $S^2$ is the unit sphere in $\R^3$. Assume that
$u=0$ on $S$, where $S$ is a bounded closed smooth connected surface in $\R^3$.

{\em Is there another bounded closed smooth connected surface of
zeros of an entire function $u_\epsilon$ of exponential type,
$u_\epsilon=\int_{S^2} e^{ik\beta \cdot x}[f(\beta)+\epsilon g(\beta)]d\beta$, where $g\in L^2(S^2)$ and $\epsilon>0$ is a small
parameter?}

 We will not use Theorem 1 since assumption \eqref{e2} may not hold, but sketch an  argument, based on the fact that $S$ in the above
question is the intersection of an analytic set with $\R^3$, see, for example, \cite{F} for the definition
and properties of analytic sets. The functions $u$ and $u_\epsilon$ in Remark 2 solve the differential equation
\be\label{e16}
\nabla^2 u +k^2 u=0 \qquad in \quad \R^3, \quad k^2=const>0.
\ee
The function $u_N$ may vanish on $S$ at most on the closed set $\sigma\subset S$ which is of the surface measure zero (by the uniqueness of the solution to the Cauchy problem for equation \eqref{e16}). For every point $s\in S\setminus \sigma$ the argument given in the proof
of Theorem 1 yields the existence of $t(s)$, the unique solution to \eqref{e6}. Since $S$ is real analytic
the set $\tilde{S}_\epsilon$, defined in the proof of Theorem 1, is analytic and $\tilde{S}_\epsilon$ is a part of the analytic set defined by the equation $u_\epsilon=0$.
In our problem $S$ is a bounded closed real analytic surface. The set $\tilde{S}_\epsilon$ can be continued analytically to
an analytic set which intersects the real space $\R^3$ over a real analytic surface $S_\epsilon$.
It is still an {\em open problem} to prove (or disprove) that the analytic continuation of the set  $\tilde{S}_\epsilon$ intersects $\R^3$ over a bounded closed real analytic surface $S_\epsilon\in \R^3$.

\end{document}